\documentclass[12pt,twoside]{article}
\usepackage{graphicx}
\usepackage{amsfonts,amsbsy,amsthm,amssymb,amsmath,tikz}
\allowdisplaybreaks
\usepackage{cite}
\usepackage{graphics}
\def\bpsp{\begin{pspicture}}
\def\epsp{\end{pspicture}}

\newtheorem{theorem}{Theorem}[section]
\newtheorem{remark}[theorem]{Remark}
\newtheorem{example}[theorem]{Example}
\newtheorem{lemma}[theorem]{Lemma}
\newtheorem{corollary}[theorem]{Corollary}
\newtheorem{definition}[theorem]{Definition}
\newtheorem{proposition}[theorem]{Proposition}

\newtheorem{note}{Note}
\newtheorem{case}{Case}
\newtheorem{conjecture}{Conjecture}
\newtheorem{question}{Question}

\newcommand{\bea}{\begin{eqnarray}}
\newcommand{\eea}{\end{eqnarray}}
\newcommand{\beq}{\begin{eqnarray*}}
\newcommand{\eeq}{\end{eqnarray*}}

\def\m4{\mbox{\rm ~(mod $4$)}}

\def \bd{\begin{definition}}
\def \ed{\end{definition}}
\def \bqu{\begin{question}}
\def \equ{\end{question}}
\def \bcc{\begin{conjecture}}
\def \ecc{\end{conjecture}}
\def \bt{\begin{theorem}}
\def \et{\end{theorem}}
\def \bl{\begin{lemma}}
\def \el{\end{lemma}}
\def \bc{\begin{corollary}}
\def \ec{\end{corollary}}
\def \be{\begin{equation}}
\def \ee{\end{equation}}
\def \ben{\begin{enumerate}}
\def \een{\end{enumerate}}
\def \ba{\begin{array}}
\def \ea{\end{array}}
\def \bp{\begin{proposition}}
\def \ep{\end{proposition}}
\def \bx{\begin{example}}
\def \ex{\end{example}}
\def \br{\begin{remark}}
\def \er{\end{remark}}
\def \bdsc{\begin{description}}
\def \edsc{\end{description}}

\def \bn{\begin{case}}
\def \en{\end{case}}
\def \bnt{\begin{note}}
\def \ent{\end{note}}
\def\1{1\!\!1}

\def\mm2{\mbox{\rm ~(mod $2$)}}
\def\m4{\mbox{\rm ~(mod $4$)}}

\def\qed{\nolinebreak\hfill\rule{.2cm}{.2cm}\par\addvspace{.5cm}}

\def\m{\mu}

\def\1{\textbf{1}}
\def\0{\textbf{0}}

\begin{document}
\title{Distance Laplacian eigenvalues of graphs and chromatic and independence number  }
\author{ S. Pirzada$ ^{a} $, Saleem Khan$ ^{b} $\\
$^{a,b}${\em Department of Mathematics, University of Kashmir, Srinagar, Kashmir, India}\\
$ ^{a} $\texttt{pirzadasd@kashmiruniversity.ac.in}; $^{b}$\texttt{khansaleem1727@gmail.com}
}
\date{}

\pagestyle{myheadings} \markboth{Pirzada, Khan}{Distance Laplacian eigenvalues of graphs}
\maketitle
\vskip 5mm
\noindent{\footnotesize \bf Abstract.} For a connected graph $G$ of order  $n$, let $Diag(Tr)$ be the diagonal matrix of vertex transmissions and $D(G)$ be the distance matrix of $G$. The distance Laplacian matrix of $G$ is defined as $D^L(G)=Diag(Tr)-D(G)$ and the eigenvalues of $D^{L}(G)$ are called the distance Laplacian eigenvalues of $G$. Let $\partial_{1}^{L}(G)\geq \partial_{2}^{L}(G)\geq \dots \geq \partial_{n}^{L}(G)$ be the distance Laplacian eigenvalues of $G$. Given an interval $I$, let $m_{D^{L} (G)} I$ (or simply $m_{D^{L} } I$) be the number of distance Laplacian eigenvalues  of $G$ which lie in the interval $I$. For a prescribed interval $I$, we determine $m_{D^{L} }I$ in terms of independence number $\alpha(G)$, chromatic number $\chi$, number of pendant vertices and diameter $d$ of the graph $G$. In particular, we prove that $m_{D^{L}(G) }[n,n+2)\leq \chi-1$, ~$m_{D^{L}(G) }[n,n+\alpha(G))\leq n-\alpha(G)$ and we show that the inequalities are sharp. We also show that $m_{D^{L} (G )}\bigg( n,n+\left\lceil\frac{n}{\chi}\right\rceil\bigg)\leq n- \left\lceil\frac{n}{\chi}\right\rceil-C_{\overline{G}}+1 $, where $C_{\overline{G}}$ is the number of components in $\overline{G}$, and discuss some cases where the bound is best possible. In addition, we prove that $m_{D^{L} (G )}[n,n+p)\leq n-p$, where $p\geq 1$ is the number of pendant vertices.  Also, we characterize graphs of diameter $d\leq 2$ which satisfy $m_{D^{L}(G) } (2n-1,2n )= \alpha(G)-1=\frac{n}{2}-1$. At the end, we propose some problems of interest.

\vskip 3mm

\noindent{\footnotesize Keywords: Distance Laplacian matrix;  distance Laplacian eigenvalues; diameter; independence number; chromatic number}

\vskip 3mm
\noindent {\footnotesize AMS subject classification: 05C50, 15A18.}

\section{Introduction}

Throughout this paper, we consider simple and connected graphs. A simple connected graph  $G=(V,E)$ consists of the vertex set $V(G)=\{v_{1},v_{2},\ldots,v_{n}\}$ and the edge set  $E(G)$. The \textit{order} and \textit{size} of $G$ are $|V(G)|=n$ and  $|E(G)|=m$, respectively. The \textit{degree} of a vertex $v,$ denoted by $d_{G}(v)$ (we simply write by $d_v$) is the number of edges incident on the vertex $v$. Further, $N_G (v)$ denotes the set of all vertices that are adjacent to $v$ in $G$ and $\overline{G}$ denotes the complement of the graph $G$. A vertex $u\in V(G)$ is called a pendant vertex if $d_{G}(u) =1$.For other standard definitions, we refer to \cite{5R8,5R9}.\\
\indent If $A$ is the adjacency matrix  and $D(G)=diag(d_1 ,d_2 ,\dots,d_n)$ is the diagonal matrix of vertex degrees of $G$, the $Laplacian$ $matrix$ of $G$ is defined as $ L(G)=D(G)-A$. By the spectrum of $G$, we mean the spectrum of its adjacency matrix, and it consists of the eigenvalues $\lambda_1 \geq \lambda_2 \geq \dots \geq \lambda_n$. The Laplacian spectrum of $G$ is the spectrum of its Laplacian matrix, and is denoted by $\mu_1 (G) \geq \mu_2 (G) \geq \dots \geq \mu_n (G) =0$. For any interval $I$, let $m_{L(G)}I$ be the number of Laplacian eigenvalues  of $G$ that lie in the interval $I$. Also, let $m_{L(G)}(\mu_i (G) )$ denote the multiplicity of the Laplacian eigenvalue $\mu_i (G) $ .\\
\indent In $G$, the \textit{distance} between the two vertices $u,v\in V(G),$ denoted by $d_{uv}$, is defined as the length of a shortest path between $u$ and $v$. The \textit{diameter} of $G$, denoted by $d$, is the maximum distance between any two vertices of $G.$ The \textit{distance matrix} of $G$, denoted by $D(G)$, is defined as $D(G)=(d_{uv})_{u,v\in V(G)}$.
The \textit{transmission} $Tr_{G}(v)$
(we will write $Tr(v)$ if the graph $G$ is understood) of a vertex $v$ is defined as the sum of the distances from $v$ to all other vertices in $G$, that is, $Tr_{G}(v)=\sum\limits_{u\in V(G)}d_{uv}.$\\
\indent Let $Tr(G)=diag (Tr(v_1),Tr(v_2),\ldots,Tr(v_n)) $ be the diagonal matrix of vertex transmissions of $G$. Aouchiche and Hansen \cite{5R1} defined the \textit{distance Laplacian matrix} of a connected graph as $D^L(G)=Tr(G)-D(G)$ (or briefly written as $D^{L}$). The eigenvalues of $D^{L}(G)$ are called the distance Laplacian eigenvalues of $G$. Since $ D^L(G) $ is a real symmetric positive semi-definite matrix, we denote its eigenvalues by $\partial_{i}^{L}(G)  $'s and order them as $0=\partial_{n}^{L}(G)\leq \partial_{n-1}^{L}(G)\leq \dots\leq \partial_{1}^{L}(G)$. The distance Laplacian eigenvalues are referred as $D^L-eigenvalues$ of $G$ whenever the graph $G$ is understood. Some recent work can be seen in \cite{pk1,pk2}. For any interval $I$,  $m_{D^L (G)}I$  represents the number of distance Laplacian eigenvalues of $G$ that lie in the interval $I$. Also, $m_{D^L  (G)}(\partial_{i}^{L}(G) )$  denotes the multiplicity of the distance Laplacian eigenvalue $ \partial_{i}^{L}(G) $.  The multiset of eigenvalues of  $ D^L(G)$ is called the \textit{distance Laplacian spectrum} of $G$. If there are only $k$ distinct distance Laplacian eigenvalues of $G$, say, $\partial_{1}^{L}(G),\partial_{2}^{L}(G),\dots,\partial_{k}^{L}(G)$ with corresponding multiplicities as $n_1 ,n_2 ,\dots, n_k$, then we convey this information in the matrix form as\\
$$\begin{pmatrix}
\partial_{1}^{L}(G) & \partial_{2}^{L}(G) & \dots & \partial_{k}^{L}(G)\\
n_1 & n_2 & \dots & n_k\\
\end{pmatrix}.$$
\indent We denote by $K_n$ the complete graph of order $n$ and by $K_{t_1 ,\dots, t_k}$ the complete multipartite graph with order of parts $t_1 ,\dots, t_k$. The star graph of order $n$ is denoted by $S_n$. Further, $SK_{n,\alpha}$ denotes the complete split graph, that is, the complement of the disjoint union of a clique $K_\alpha$ and $n-\alpha$ isolated vertices.
For two disjoint graphs $G$ and $H$ of order $n_1$ and $n_2$, respectively, the \textit{corona graph} $GoH$ is the graph obtained by taking one copy of $G$ and $n_1$ copies of $H$, and then joining the \textit{i}th vertex of $G$ to every vertex in the \textit{i}th copy of $H$, for all $ 1\leq i\leq n_1$.\\
\indent In a graph $G$, the subset $M\subseteq V(G)$ is called an \textit{independent set} if no two vertices of $M$ are adjacent. The \textit{independence number} of $G$ is the cardinality of the largest independent set of $G$ and is denoted by $\alpha(G)$. A set $M\subseteq V(G)$ is \textit{dominating} if every $v\in V(G) \setminus M$ is adjacent to some member in $S$. The \textit{domination number} $\gamma(G)$ is the minimum size of a dominating set.\\
\indent The \textit{chromatic number} of a graph $G$ is the minimum number of colors required to color the vertices of $G$ such that no two adjacent vertices get the same color. It is denoted by $\chi(G)$. The set of all vertices with the same color is called a \textit{color class}. \\
\indent The distribution of Laplacian eigenvalues of a graph $G$ in relation to various graph parameters of $G$ has been studied extensively. Grone and Merris  \cite{5R10} and Merris \cite{5R11} obtained bounds for $m_{L(G)}[0,1) $ and $m_{L(G)}[0,2) $. Guo and Wang \cite{5R12} showed that if $G$ is a connected graph with matching number $\nu(G)$, then $m_{L(G)}(2,n]>\nu(G)$, where $n>2\nu(G)$. Some work in this direction can be seen in \cite{cjt}. Recently, Ahanjideh et al \cite{5R0} obtained bounds for $m_{L(G)}I $ in terms of structural parameters of $G$. In particular, they showed that $m_{L(G)}(n -\alpha(G), n] \leq n -\alpha(G)$ and $m_{L(G)}(n-d(G)+3, n]\leq n -d(G) -1$, where $\alpha(G)$ and $d(G)$ denote the independence number and the diameter of $G$, respectively. The distribution of the distance Laplacian eigenvalues of a graph $G$ with respect to its structural parameters has not got its due attention and our investigation in this manuscript  is an attempt in that direction.  \\
\indent The rest of the paper is organized as follows.  In Section 2, we find the distribution of $G$ in relation to the chromatic number $\chi$ and the number of pendant vertices. We show that $m_{D^{L}(G) }[n,n+2)\leq \chi-1$ and show that the inequality is sharp. We also prove that  $m_{D^{L} (G )}\bigg( n,n+\left\lceil\frac{n}{\chi}\right\rceil\bigg)\leq n- \left\lceil\frac{n}{\chi}\right\rceil-C_{\overline{G}}+1 $, where $C_{\overline{G}}$ is the number of components in $\overline{G}$, and discuss some cases where the bound is best possible. In addition, we prove that $m_{D^{L} (G )}[n,n+p)\leq n-p$, where $p\geq 1$ is the number of pendant vertices. In Section 3, we determine the distribution of distance Laplacian eigenvalues  of $G$ in terms of the independence number $\alpha(G)$ and diameter $d$. In particular, we show that   $m_{D^{L} (G)}[n,n+\alpha(G))\leq n-\alpha(G)$ and show that the inequality is sharp. We show that $m_{D^{L}(G)}[0,dn]\geq d+1$.   We characterize the graphs having diameter $d\leq 2$ satisfying  $m_{D^{L}(G) } (2n-1,2n )= \alpha(G)-1=\frac{n}{2}-1$. In Section 4, we propose some research problems.

\section{Distribution of distance Laplacian eigenvalues, chromatic number and pendant vertices }

For a graph $G$ with $n$ vertices, let $Tr_{max}(G)=max\{Tr(v):v\in V(G)\}$ . Whenever the graph $G$ is understood, we will write $Tr_{max}$ in place of $Tr_{max}(G)$. We have the following important result from matrix theory.

\begin{lemma}\label{L2}\emph {\cite{5R3}} Let $M=(m_{ij})$ be a $n\times n$ complex matrix having $l_1 ,l_2 ,\dots,l_p$ as its distinct eigenvalues. Then
$$\{l_1 ,l_2 ,\dots,l_p\}\subset \bigcup\limits_{i=1}^{n}\Big \{z:|z-m_{ii}|\leq \sum\limits_{j\neq i}|m_{ij}|\Big\}.$$

\end{lemma}
\indent By using Lemma \ref{L2} for the distance Laplacian matrix of a graph $G$ with $n$ vertices, we get
\begin{equation}
\partial^L_{1}(G)\leq 2Tr_{max}
\end{equation}
The following fact about distance Laplacian eigenvalues will be used in the sequel.\\
\textbf{Fact 1.} Let $G$ be a connected graph of order $n$ and having distance Laplacian eigenvalues in the order $\partial^L_{1}(G)\geq \partial^L_{2}(G)\geq \dots \geq \partial^L_{n}(G)$. Then,\\
\hspace*{25mm} $\partial^L_{n}(G)=0$ and $\partial^L_{i}(G)\geq n$ for all $i=1,2,\dots,n-1.$\\\\
 We recall  the following important results.

 \begin{theorem}\label{T7}
(Cauchy Interlacing Theorem). Let $M$ be a real symmetric matrix of order $n$, and let $A$ be a principal submatrix of $M$ with order $s\leq n$. Then $$\lambda_i (M)\geq \lambda_i (A) \geq \lambda_{i+n-s} (M)\hspace{1cm}(1\leq i\leq s).$$
\end{theorem}
\begin{lemma} \label{L1}\emph {\cite{5R1}}  Let $G$ be a connected graph with $n$ vertices and $m$ edges, where $m\geq n$. Let $G^*$ be the connected graph obtained from $G$ by deleting an edge. Let $\partial^L_1 \geq \partial^L_2 \geq ...\geq \partial^L_n$ and ${\partial^*_1}^L \geq {\partial^*_2}^L \geq ...\geq {\partial^*_n}^L$ be the spectrum of $G$ and $G^*$, respectively. Then ${\partial^*_i}^L \geq \partial^L_i $ for all $i=1,\dots,n$.
\end{lemma}
\begin{lemma}\label{L8} \emph{\cite{5R7} } Let $t_{1},t_{2},\dots,t_{k}$ and n be integers such that $t_{1}+t_{2}+\dots+t_{k}=n$ and $t_{i}\geq 1$ for $i=1,2,\dots,k$. Let $p=|\{i:t_{i}\geq 2\}|$. The distance Laplacian spectrum of the complete $k-partite$ graph $K_{t_{1},t_{2},\dots,t_{k}}$ is$ \Big((n+t_{1})^{(t_{1}-1)},\dots,(n+t_{p})^{(t_{p}-1)},n^{(k-1)},0\Big)$.
\end{lemma}
\begin{lemma}\label{L3} \emph {\cite{5R1}} Let $G$ be a connected graph with $n$ vertices. Then $\partial^L_{n-1}\geq n$ with equality if and only if $\overline{G}$ is disconnected. Furthermore, the multiplicity of $n$ as an eigenvalue of $D^L (G)$ is one less than the number of components of  $\overline{G}$.
\end{lemma}

First we obtain an upper bound for $m_{D^{L} (G)} I$, where $I$ is the interval $[n,n+2)$, in terms of the chromatic number $\chi$ of $G$.

\begin{theorem} \label{T8} Let $G$ be a connected graph of order $n$ and having chromatic number $\chi$. Then $$m_{D^{L} (G)} [n,n+2 ) \leq \chi-1.$$ Inequality is sharp and is shown by all complete multipartite graphs.
\end{theorem}
 \noindent {\bf Proof.}  Let $t_1 ,t_2 ,\dots,t_\chi $ be $\chi$  positive integers such that $t_1 +t_2 +\dots+t_{\chi} =n$ and let these numbers be the cardinalities of $\chi$ partite classes of $G$. We order these numbers as $t_1 \geq t_2 \geq \dots\geq t_{\chi(G)} $. Thus $G$ can be considered as a spanning subgraph of the complete multipartite graph $H=K_{t_1 ,t_2 ,\dots,t_{\chi}}$ with  $t_1 \geq t_2 \geq \dots\geq t_{\chi} $ as the cardinalities of its partite classes. Using Lemma \ref{L8}, we see that $m_{D^{L} (H )} [n,n+2 ) = \chi-1$. By Lemma \ref{L1} and the Fact 1, we have  $ m_{D^{L} (G )} [n,n+2 ) \leq m_{D^{L} (H )} [n,n+2 ) = \chi-1$, proving the inequality. Using Lemma \ref{L8}, we see that the equality holds for all complete multipartite graphs. \qed

 As a consequence of Theorem \ref{T8}, we have the following observation.

 \begin{corollary} \label{C2}   Let $G$ be a connected graph of order $n$  having chromatic number $\chi$. Then $$ m_{D^{L} (G )} [n+2,2Tr_{max} ]\geq n- \chi.$$ Inequality is sharp and is shown by all complete multipartite graphs.
 \end{corollary}
 \noindent {\bf Proof.} By using the Fact 1, we get
 \begin{align*}
  &m_{D^{L} (G )} [n,n+2 )+ m_{D^{L} (G )}[n+2,2Tr_{max} ]  =n-1, \\&
or ~~~~~  \chi-1+ m_{D^{L} (G )}[n+2,2Tr_{max} ] \geq n-1, \\&
or ~~~~~~   m_{D^{L} (G )}[n+2,2Tr_{max} ] \geq  n- \chi.
\end{align*}
Therefore, the inequality is established. The remaining part of the proof follows from Theorem \ref{T8}. \qed

In the following theorem, we characterize the unique graph with chromatic classes of the same cardinality having $n-1$ eigenvalues in the interval $\big[n,n+\frac{n}{\chi}\big]$.

\begin{theorem} \label{T9}  Let $G$ be a connected graph of order $n$ and having the chromatic number $\chi$. If the chromatic classes are of the same cardinality, then
 $$ m_{D^{L} (G )} \big[n,n+\frac{n}{\chi}\big]\leq n-1$$ with equality if and only if $G\cong K_{\frac{n}{\chi},\dots,\frac{n}{\chi}}$.
\end{theorem}
 \noindent {\bf Proof.} Using Fact 1, we get the required inequality. Now, we will show that the equality holds for the graph $H= K_{\frac{n}{\chi},\dots,\frac{n}{\chi}}$. Using Lemma \ref{L8}, we have the distance Laplacian spectrum of $H$ as
 $$\begin{pmatrix}
0 & n & n+\frac{n}{\chi} \\
1 & \chi-1 & n-\chi \\
\end{pmatrix},$$
which clearly shows that the equality holds for the graph $H$. To complete the proof, we will show that if $G\ncong H$, then $ m_{D^{L} (G )} \big[n,n+\frac{n}{\chi}\big]< n-1$. Since the chromatic classes are of the same cardinality, we see that $G$ has to be an induced subgraph of $H$ and $n=s\chi$ for some integer $s$, so that $s=\frac{n}{\chi}$. In $H$, let $e=\{u,v\}$ be an edge  between the vertices $u $ and $v$. Using Lemma \ref{L1}, it is sufficient to take $G=H-e$. In $G$, we see that $Tr(u)=Tr(v)=n+s-1$. Let $A$ be the principal submatrix of $D^L (G)$  corresponding to the vertices $u$ and $v$. Then $A$ is given by
\begin{equation*}
A=
\begin{bmatrix}
n+s-1 & -2 \\
-2 & n+s-1
\end{bmatrix}.
\end{equation*}
Let $c(x)$ be the characteristic polynomial of $A$. Then $c(x)=x^2 -2(n+s-1)x+{(n+s-1)}^2-4$. Let $x_1 $ and $x_2$ be the roots of $c(x)$ with $x_1 \geq x_2$. It can be easily seen that $x_1=n+s+1$. Using Theorem \ref{T7}, we have $\partial^L _1 (G)\geq x_1 =n+s+1>n+s=n+\frac{n}{\chi}$. Thus,  $ m_{D^{L} (G )} \big[n,n+\frac{n}{\chi}\big]< n-1$ and the proof is complete. \qed

Now, we obtain an upper bound for the number of distance Laplacian eigenvalues which fall in the interval $\bigg( n,n+\left\lceil\frac{n}{\chi}\right\rceil\bigg)$.

\begin{theorem}\label{TN1}Let $G\ncong K_n$ be a connected graph on $n$ vertices with chromatic number $\chi$. Then,
\begin{equation}
m_{D^{L} (G )}\bigg( n,n+\left\lceil\frac{n}{\chi}\right\rceil\bigg)\leq n- \left\lceil\frac{n}{\chi}\right\rceil-C_{\overline{G}}+1
\end{equation}
where $C_{\overline{G}}$ is the number of components in $\overline{G}$. The bound is best possible for  $\chi=2$ (when $n$ is odd) and $\chi=n-1$   as shown by  $K_{m+1, m}$, where $n=2m+1$, and $K_{2,\underbrace{1,1,\dots,1}_{n-2}} $, respectively.
\end{theorem}
\noindent {\bf Proof.}  Let $n_1 \geq n_2 \geq \dots\geq n_{\chi} $ be $\chi$  positive integers in that order such that $n_1 +n_2 +\dots+n_{\chi} =n$ and let these numbers be the cardinalities of $\chi$ partite classes of $G$. Clearly, $G$ can be considered as a spanning subgraph of the complete multipartite graph $H=K_{n_1 ,n_2 ,\dots,n_{\chi}}$. Using Lemmas \ref{L1} and \ref{L8}, we get
$$\partial^L _i (G)\geq \partial^L _i (H)=n+n_1, ~~~~~~ \text{for all}  ~ 1\leq i\leq n_1 -1.$$
 As $n_1$ is largest among the  cardinalities of chromatic classes, it is at least equal to average, that is,
 $n_1 \geq \frac{n}{\chi}$. Also, $n_1$ is an integer, therefore, $n_1 \geq \left\lceil\frac{n}{\chi}\right\rceil$. Using this fact in above inequality, we get
$$
 \partial^L _i (G)\geq n+\left\lceil\frac{n}{\chi}\right\rceil ~~~~~~ \text{ for all} ~ 1\leq i\leq n_1 -1.
 $$
 Thus, there at least $n_1 -1$ distance Laplacian eigenvalues of $G$ which are greater than or equal to $n+\left\lceil\frac{n}{\chi}\right\rceil$.
Also from Lemma \ref{L3}, we see that $n$ is a distance Laplacian eigenvalues of $G$ with multiplicity exactly $C_{\overline{G}}-1$. Using these observations with Fact 1, we get
\begin{align*}
m_{D^{L} (G )}\bigg( n,n+\left\lceil\frac{n}{\chi}\right\rceil\bigg)& \leq n- (n_1 -1)-(C_{\overline{G}}-1)-1\\
& = n-n_1 -C_{\overline{G}}+1\\
& \leq n-\left\lceil\frac{n}{\chi}\right\rceil-C_{\overline{G}}+1,
\end{align*}
proving the required inequality. \\
 Let $G^*=K_{2,\underbrace{1,1,\dots,1}_{n-2}} $. It is easy to see that $\left\lceil\frac{n}{n-1}\right\rceil=2$. Also, the complement of $G^*$ has exactly $n-1$ components. By Lemma \ref{L8},  the distance Laplacian spectrum of $G^*$ is given as follows
 $$\begin{pmatrix}
0 & n & n+2 \\
1 & n-2 & 1 \\
\end{pmatrix}.$$
Putting all these observations in Inequality (2.2), we see that the equality holds for $G^*$ which shows that the bound is best possible when $\chi=n-1$.

Let $G^{**}=K_{m+1, m}$, where $n=2m+1$. In this case, we see that $\left\lceil\frac{n}{2}\right\rceil=m+1=\frac{n+1}{2}$ and the complement of $G^{**}$ has exactly $2$ components.  By Lemma \ref{L8}, we observe that the distance Laplacian spectrum of $G^{**}$ is given as follows
 $$\begin{pmatrix}
0 & n & \frac{3n+1}{2} & \frac{3n-1}{2} \\
1 & 1 & \frac{n-1}{2} & \frac{n-3}{2}\\
\end{pmatrix}.$$
 Using all the above observations in Inequality (2.2), we see that the equality holds for $G^{**}=K_{m+1, m}$ which shows that the bound is best possible when $\chi=2$ and $n$ is odd. \qed

 The following are some immediate consequences of Theorem \ref{TN1}.

\begin{corollary}\label{CN1} Let $G\ncong K_n$ be a connected graph on $n$ vertices with chromatic number $\chi$. Then,
$$
m_{D^{L} (G )}\bigg[ n+\left\lceil\frac{n}{\chi}\right\rceil,\partial^L _1 (G)\bigg]\geq \left\lceil\frac{n}{\chi}\right\rceil-1.
$$
The bound is best possible for  $\chi=2$ (when $n$ is odd) and $\chi=n-1$   as shown by  $K_{m+1, m}$, where $n=2m+1$, and $K_{2,\underbrace{1,1,\dots,1}_{n-2}} $, respectively.
\end{corollary}

\begin{corollary}\label{CN2}Let $G\ncong K_n$ be a connected graph on $n$ vertices with chromatic number $\chi$. If $\overline{G}$ is connected, then
$$m_{D^{L} (G )}\bigg( n,n+\left\lceil\frac{n}{\chi}\right\rceil\bigg)\leq n- \left\lceil\frac{n}{\chi}\right\rceil.$$
\end{corollary}
\noindent{\bf Proof.} Since $\overline{G}$ is connected, therefore, $C_{\overline{G}}=1$. Putting $C_{\overline{G}}=1$ in Inequality (2.2) proves the desired result. \qed

The next theorem shows that there are at most $n-p$ distance Laplacian eigenvalues of $G$ in the interval $[n,n+p)$, where $p\geq 1$ is the number of pendant vertices in $G$.

\begin{theorem}\label{TN2} Let $G\ncong K_n$ be a connected graph on $n$
vertices having $p\geq 1$ pendant vertices, then
$$m_{D^{L} (G )}[n,n+p)\leq n-p.$$
For $p=n-1$, equality holds if and only if $G\cong S_n$.
\end{theorem}
\noindent{\bf Proof.} Let $S$ be the set of pendant vertices so that $|S|=p$. Clearly, $S$ is an independent set of $G$.  Obviously, the induced subgraph, say $H$, on the vertex set $M=V(G)\setminus S$ is  connected. Let  the chromatic number of $H$ be $q$ and  $n_1 \geq n_2 \geq \dots \geq n_q$ be the cardinalities of these chromatic classes in that order, where $1\leq q \leq n-p$ and $n_1 +n_2 +\dots+n_q =n-p$. Let $n_k \geq p \geq n_{k+1}$, where $0\leq k \leq q$, $n_0 =p$ if  $k=0$ and $n_{q+1}=p$ if $k=q$. With this partition of the vertex set $V(G)$ into $q+1$ independent sets, we easily see that $G$ can be considered as an induced subgraph of complete $q+1$-partite graph $L=K_{n_1 ,n_2,\dots, n_k ,p,n_{k+1} ,\dots,n_q} $. Consider the following two cases.\\
\noindent{\bf Case 1.} Let $1\leq k \leq q$ so that $n_1 \geq p$. Then, from Lemmas \ref{L1} and \ref{L8}, we get
$$\partial^L _i (G)\geq \partial^L _i (L)=n+n_1\geq n+p, ~~~ \text{ for all} ~  1\leq i \leq n_1 -1. $$
\noindent{\bf Case 2.} Let $k=0$ so that $p\geq n_1$. Again, using Lemmas \ref{L1} and \ref{L8}, we get
$$\partial^L _i (G)\geq \partial^L _i (L)=n+p, ~~~ \text{ for all} ~  1\leq i \leq p -1.$$
Thus, in both cases, we see that there are at least $p-1$ distance Laplacian eigenvalues of $G$ which are greater than or equal to $n+p$. As $p\geq 1$, so $\overline{G}$ has at most two components, which after using Lemma \ref{L3} shows that $n$ is a distance Laplacian eigenvalue of $G$ of multiplicity at most one. From the above observations and Fact 1, we get
 $$m_{D^{L} (G )}[n,n+p)\leq n-p,$$
 which proves the required inequality.

For the second part of the theorem, we see that $S_n$ is the only connected graph having $n-1$ pendant vertices. The distance Laplacian spectrum of $S_n$ by Lemma \ref{L8} is given as
  $$\begin{pmatrix}
0 & n & 2n-1 \\
1 & 1 & n-2\\
\end{pmatrix}$$
and the proof is complete. \qed

An immediate consequence is as follows.

\begin{corollary}\label{CN3}
 Let $G\ncong K_n$ be a connected graph on $n$
vertices having $p\geq 1$ pendant vertices, then
$$m_{D^{L} (G )}[n+p,\partial^L _1 (G)]\geq p-1.$$
For $p=n-1$, equality holds if and only if $G\cong S_n$.
\end{corollary}
The following lemma will be used in the proof of Theorem \ref{T11}.

\begin{lemma}\label{L9}  \emph{\cite{5R2}} Let $G$ be a graph with $n$ vertices. If $K=\{v_1 ,v_2 ,\dots,v_p\}$ is an independent set of $G$ such that $N(v_i)=N(v_j)$ for all $i,j\in \{1,2,\dots,p\}$, then $\partial=Tr(v_i)=Tr(v_j)$ for all $i,j\in \{1,2,\dots,p\}$ and $\partial +2$ is an eigenvalue of $D^L (G)$ with multiplicity at least $p-1$.
\end{lemma}

\begin{theorem} \label{T11} Let $G$ be a connected graph of order $n\geq 4$ having chromatic number $\chi$. If $S=\{v_1 ,v_2 ,\dots,v_p\} \subseteq V(G)$, where $|S|=p\geq \frac{n}{2}$, is the set of pendant vertices such that every vertex in $S$ has the same neighbour in $V(G)\setminus S$, then
$$ m_{D^{L} (G )} [n,2n-1)\leq n-\chi.$$
\end{theorem}
\noindent {\bf Proof.} Clearly, all the vertices in $S$ form an independent set. Since all the vertices in $S$ are adjacent to the same vertex, therefore, all the vertices of $S$ have the same transmission. Now, for any  $v_i$  $(i=1,2,\dots,p)$ of $S$, we have
\begin{align*}
T=Tr(v_i ) \geq 2(p-1)+1+2(n-p-1) =2n-3.
\end{align*}
From Lemma \ref{L9}, there are at least $p-1$ distance Laplacian eigenvalues of $G$ which are greater than or equal to $T+2$. From above, we have $T+2\geq 2n-3+2=2n-1$. Thus,  there are at least $p-1$ distance Laplacian eigenvalues of $G$ which are greater than or equal to $2n-1$, that is, $ m_{D^{L} (G )} [2n-1,2Tr_{max}]\geq p-1$. Using Fact 1, we have
\begin{equation}
  m_{D^{L} (G )} [n,2n-1)\leq n-p.
  \end{equation}
  We  claim that $\chi(G)\leq \frac{n}{2}$. If possible, let $\chi(G)> \frac{n}{2}$. We have following two cases to consider.\\
  $\bf {Case ~ 1.}$ Let $p=n-1$. Clearly, the star is the only connected graph having $n-1$ pendant vertices. Thus, $G\cong S_n$. Also, $\chi(S_n)=2$, a contradiction, as $\chi(S_n)=2\leq\frac{n}{2}$, for $n\geq 4$.\\
 $\bf {Case ~ 2.}$ $\frac{n}{2}\leq p \leq n-2$. Since $p\leq n-2$, there is at least one vertex, say $u$, which is not adjacent to any vertex in $S$. Thus in the minimal coloring of $G$, at least $p+1$ vertices, say, $u,v_1 ,\dots,v_p$ can be colored using only one color. The rest $n-p-1$ vertices can be colored with at most $n-p-1$ colors. Thus, $\chi\leq 1+n-p-1=n-p\leq n-\frac{n}{2}=\frac{n}{2}$, a contradiction. Therefore, $\chi \leq \frac{n}{2}\leq p$. Using this in Inequality (2.3), we get
 $$ m_{D^{L} (G )} [n,2n-3)\leq n-\chi,$$
completing the proof.  \qed

 To have a bound only  in terms of order $n$ and the number of pendant vertices $p$, we can relax the conditions $p\geq \frac{n}{2}$ and $n\geq 4$ in  Theorem \ref{T11}. This is given in the following corollary.

 \begin{corollary} \label{C3} Let $G$ be a connected graph of order $n$ . If $S=\{v_1 ,v_2 ,\dots,v_p\} \subseteq V(G)$ is the set of pendant vertices such that  every vertex in $S$ has the same neighbour in $V(G)\setminus S$, then
$$ m_{D^{L} (G )} [n,2n-1)\leq n-p.$$
 \end{corollary}

\section{Distribution of distance Laplacian eigenvalues, independence number and diameter}

The following lemma will be useful.

Now, we obtain an upper bound for $m_{D^{L} (G)}I$, where $I$ is the interval $[n,n+\alpha(G))$, in terms of order $n$ and independence number $\alpha(G)$.

\begin{theorem} \label{T1} Let $G$ be a connected graph of order $n$ having independence number $\alpha (G)$. Then $m_{D^{L} (G)} [n,n+\alpha(G))\leq n-\alpha(G)$. For $\alpha(G)=1$ or $\alpha(G)=n-1$,  the  equality holds if and only if $G\cong K_n$ or $G\cong S_n$. Moreover,  for every integer $n$ and $\alpha(G)$ with $2\leq \alpha(G)\leq n-2$, the bound is sharp, as $SK_{n,\alpha}$ satisfies the inequality.
\end{theorem}
\noindent {\bf Proof.} We have the following three cases to consider.\\
{\bf Case 1.}  $\alpha(G)=1$. Clearly, in this case $G\cong K_n$ and the distance Laplacian spectrum of a complete graph is
$$\begin{pmatrix}
0 & n \\
1 & n-1 \\
\end{pmatrix}.$$
Therefore, we have $m_{D^{L} (K_n)} [n,n+1)= n-1$ which proves the result in this case. \\
 {\bf Case 2.}  $\alpha(G)= n-1$. Since the star $S_n$ is the only connected graph having independence number $n-1$, therefore, $G\cong S_n$ in this case. Now, $n-\alpha(S_n)=n-n+1=1$. From Lemma \ref{L8}, the distance Laplacian spectrum of $S_n $ is given as  \\
  $$\begin{pmatrix}
0 & n & 2n-1 \\
1 & 1 & n-2 \\
\end{pmatrix}.$$
Therefore, $m_{D^{L} (S_n)} [n,2n-1)=  1$, proving the result in this case.\\
{\bf Case 3.}  $2\leq \alpha(G)\leq n-2$.  Without loss of generality, assume that $N=\{v_1 ,v_2 ,\dots ,v_{\alpha(G)}\} \subseteq V(G)$ is an independent set with maximum cardinality.  Let $H$ be the new graph obtained by adding edges between all non-adjacent vertices  in $V(G)\setminus N$ and adding edges between each vertex of $N$ to vertex of $V(G)\setminus N$. With this construction, we see that $H\cong SK_{n,\alpha}$.   Using Fact 1 and  Lemma \ref{L1}, we see that $m_{D^{L} (G)} [n,n+\alpha(G))\leq m_{D^{L} (H)} [n,n+\alpha(G))$. So to complete the proof in this case, it is sufficient to prove that $ m_{D^{L} (H)} [n,n+\alpha(G))\leq n-\alpha(G)$. By Corollary 2.4 in \cite{5R2}, the distance Laplacian spectrum of $H$ is given by
 $$\begin{pmatrix}
0 & n & n+\alpha(G) \\
1 & n-\alpha(G) & \alpha(G)-1  \\
\end{pmatrix}.$$
This shows that  $ m_{D^{L} (H)} [n,n+\alpha(G))= n-\alpha(G)$. Thus the bound is established. Also, it is clear that $SK_{n,\alpha}$ satisfies the inequality for $2\leq \alpha(G)\leq n-2$. \qed

From Theorem \ref{T1}, we have the following observation.

\begin{corollary} \label{c1} If $G$ is a connected graph of order $n$ having independence number $\alpha (G)$, then $\alpha(G) \leq 1+m_{D^{L} (G)} [n+\alpha(G),2Tr_{max}]$. For $\alpha(G)=1$ or $\alpha(G)=n-1$, the equality holds if and only if $G\cong K_n$ or $G\cong S_n$. Moreover,  for every integer $n$ and $\alpha(G)$ with $2\leq \alpha(G)\leq n-2$, the bound is sharp, as $SK_{n,\alpha}$ satisfies the inequality.
\end{corollary}
\noindent {\bf Proof.} Using Inequality (2.1) and Theorem \ref{T1}, we have
 \begin{align*}
  & m_{D^{L} (G)} [n,n+\alpha(G))+ m_{D^{L} (G)} [n+\alpha(G),2Tr_{max}]=n-1\\
 or ~~~~~~~~~~~~~~ & ~ n-\alpha(G)+  m_{D^{L} (G)} [n+\alpha(G),2Tr_{max}]\geq n-1\\
 or ~~~~~~~~~~~~~~ & ~  \alpha(G) \leq 1+m_{D^{L} (G)} [n+\alpha(G),2Tr_{max}],
 \end{align*}
 which proves the inequality. The proof of the remaining part is similar to the proof of Theorem \ref{T1}. \qed

The next result is an upper bound for $ m_{D^{L} (G)} (n,n+\alpha(G))$ in terms of the independence number $\alpha(G)$, order $n$ and number of components of the complement $\overline{G}$ of $G$.

\begin{theorem} \label{T2} Let $G$ be a connected graph with $n$ vertices having independence number $\alpha(G)$.  Then
$$ m_{D^{L} (G)} (n,n+\alpha(G))\leq n-\alpha(G) +1-k,$$
where $k$ is the number of components of $\overline{G}$. For $\alpha(G)=1$ or $\alpha(G)=n-1$, equality holds if and only if $G\cong K_n$ or $G\cong S_n$. Furthermore,  for every integer $n$ and $\alpha(G)$ with $2\leq \alpha(G)\leq n-2$, the bound is sharp, as $SK_{n,\alpha}$ satisfies the inequality.
\end{theorem}
\noindent {\bf Proof.} Since  $\overline{G}$ has $k$ components, therefore by Lemma \ref{L3}, $n$ is a distance Laplacian eigenvalue of multiplicity exactly $k-1$. Using Theorem \ref{T1}, we have
\begin{align*}
m_{D^{L} (G)} (n,n+\alpha(G)) & =m_{D^{L} (G)} [n,n+\alpha(G))-m_{D^{L} (G)} (n)\\
& =m_{D^{L} (G)} [n,n+\alpha(G))-k+1\\
& \leq n-\alpha(G) +1-k.
\end{align*}
Thus the inequality is established. The remaining part of the proof follows by observing the distance Laplacian spectrum of the graphs $ K_n$,  $ S_n$ and $SK_{n,\alpha}$ given in Theorem \ref{T1}.  \qed

We will use the following lemmas in the proof of Theorem \ref{T3}.

\begin{lemma} \label{L4} \emph {\cite{5R4}} If $G$ is a graph with domination number $\gamma (G)$, then $ m_{L(G)} [0,1)\leq \gamma (G) $.
\end{lemma}

\begin{lemma}\label{L5}\emph{\cite{5R1}} Let $G$ be a connected graph with $n$ vertices and diameter $d(G)\leq 2$. Let $\mu_1 (G) \geq \mu_2 (G)\geq \dots \geq \mu_n (G)=0$ be the Laplacian spectrum of $G$. Then the distance Laplacian spectrum of $G$ is  $2n-\mu_{n-1} (G) \geq 2n- \mu_{n-2} (G)\geq \dots \geq 2n-\mu_1 (G)>\partial^L_n (G)=0$. Moreover, for every $i\in \{1,2,\dots,n-1\}$, the eigenspaces corresponding to $\mu_i (G)$ and $2n-\mu_i (G)$ are same.
\end{lemma}

Now, we obtain an upper bound for $m_{D^{L}(G) }$, where $I$ is the interval $(2n-1,2n)$, in terms of the independence number $\alpha(G)$. This upper bound is for graphs with diameter $d(G)\leq 2$.

\begin{theorem} \label{T3} Let $G$ be a connected graph with $n$ vertices having independence number $\alpha(G)$ and diameter  $d(G)\leq 2$. Then
$$m_{D^{L} (G)} (2n-1,2n )\leq \alpha(G) -1$$and inequality is sharp as shown by $K_n$.
\end{theorem}
\noindent {\bf Proof.} We know that every maximal independent set of a graph $G$ is a minimal dominating set of $G$. Therefore, $\alpha(G)\leq \gamma (G)$. Using Lemma \ref{L4}, we get $\alpha(G)\geq  m_{L(G)} [0,1)$. As $G$ connected, the multiplicity of 0 as a Laplacian eigenvalue of $G$ is one. Thus, $\alpha(G)-1\geq  m_{L(G)} (0,1)$, that is, there are at least $\alpha(G)-1$ Laplacian eigenvalues of $G$ which are greater than zero and less than one. Using this fact in Lemma \ref{L5}, we observe that there are at least $\alpha(G)-1$ distance Laplacian eigenvalues of $G$ which are greater than $2n-1$ and less than $2n$. Thus,
$$m_{D^{L} (G)} (2n-1,2n )\leq \alpha(G) -1.$$
Clearly, $ m_{D^{L} (K_n)} (2n-1,2n )=0$ and $\alpha(K_n)=1$, which shows that equality holds for $K_n$. \qed

Our next result shows that the  upper bound in Theorem \ref{T3} can be improved for the graphs having independence number greater than $\frac{n}{2}$.

\begin{theorem} \label{L7}  Let $G$ be a connected graph with $n$ vertices having independence number $\alpha(G)>\frac{n}{2}$ and diameter  $d(G)\leq 2$ . Then $m_{D^{L} (G)} (2n-1,2n )\leq \alpha(G) -2.$
\end{theorem}
\noindent {\bf Proof.} If possible, let $m_{D^{L} (G)} (2n-1,2n )\geq \alpha(G) -1$. Using Lemma \ref{L5}, we see that there are at least $\alpha(G) -1$ Laplacian eigenvalues of $G$ which are greater than zero and less than one. As $G$ is connected, 0 is a Laplacian eigenvalue of multiplicity one. Using these facts and Lemma \ref{L4}, we have $\alpha(G) \leq m_{L(G)} [0,1)\leq \gamma(G) \leq \alpha(G).$ Thus, $ \gamma(G) =\alpha(G) >\frac{n}{2}$. This contradicts the well known fact that $\gamma(G) \leq \frac{n}{2}$ . Thus the result is established.\qed

We also use the following lemma in our next result.

\begin{lemma} \label{L6} \emph{\cite{5R5}} Let $G$ and $G^*$ be graphs with $n_1$ and $n_2$ vertices, respectively. Assume that $\mu_1 \leq \dots \leq \mu_{n_1 }$ and $\lambda_1 \leq \dots \leq \lambda_{n_2 }$ are the Laplacian eigenvalues of $G$ and $G^*$ , respectively. Then the Laplacian spectrum of $GoG^*$ is given as follows.\\
(i) The eigenvalue $\lambda_j +1$ with multiplicity $n_1$ for every eigenvalue $\lambda_j (j=2,\dots,n_2)$ of $G^*$;\\
(ii) Two multiplicity-one eigenvalues $\frac{\mu_i +n_2 +1\pm \sqrt{{(\mu_i +n_2 +1)}^2-4\mu_i}}{2}$, for each eigenvalue $\mu_i (i=1,\dots,n_1)$ of $G$.
\end{lemma}

The following result characterizes the graphs with diameter $d(G)\leq 2$ and independence number $\alpha(G)$ which satisfy $m_{D^{L} } (2n-1,2n )= \alpha(G)-1=\frac{n}{2}-1$.

\begin{theorem} \label{T5}  Let $G$ be a connected graph with $n$ vertices having independence number $\alpha(G)$ and diameter $d(G)\leq 2$. Then $m_{D^{L} (G)} (2n-1,2n )= \alpha(G)-1=\frac{n}{2}-1$ if and only if $G=HoK_1$ for some connected graph $H$.
\end{theorem}
\noindent {\bf Proof.} Assume that $G=HoK_1$ for some connected graph $H$. Then $|H|=\frac{n}{2}$. Let the Laplacian eigenvalues of $H$ be  $\mu_1 \geq \dots \geq \mu_{\frac{n}{2}}$. By Lemma \ref{L6}, the Laplacian eigenvalues of $G$ are equal to $\frac{\mu_i +2\pm \sqrt{{\mu_i}^2 +4}}{2} $, $i=1,\dots,\frac{n}{2}$. We observe that half of these eigenvalues are greater than 1 and the other half are less than 1. As $G$ is connected, 0 is a Laplacian eigenvalue of multiplicity one. So  $m_{{L} (G)} (0,1 )=\frac{n}{2}-1$. Using Lemma \ref{L5}, we see that there are $\frac{n}{2}-1$ distance Laplacian eigenvalues which are greater than $2n-1$ and less than $2n$. Thus,  $m_{D^{L} (G)} (2n-1,2n )= \frac{n}{2}-1$. Now, we will show that $\alpha(G)=\frac{n}{2}  $.  Assume that $V(G)=\{v_1, \dots,v_{\frac{n}{2}}, v'_1 ,\dots,v'_{\frac{n}{2}}\},$ where $V(H)=\{v_1, \dots,v_{\frac{n}{2}}\}$ and $N_G (v'_i)=\{v_i \}$. If $A$ is a maximal independent set, then $|A|\leq \frac{n}{2}$. For if $|A|> \frac{n}{2}$, then from the structure of $G$, we have at least one pair of vertices in $A$, say $v_i ,v'_i$, which are adjacent, a contradiction. As $\{ v'_1 ,\dots,v'_{\frac{n}{2}}\}$ is an independent set, therefore $\alpha(G)=\frac{n}{2}$. Thus, we have $m_{D^{L} (G)} (2n-1,2n )= \alpha(G)-1=\frac{n}{2}-1$.\\
\indent Conversely, assume  that $m_{D^{L} (G)} (2n-1,2n )= \alpha(G)-1=\frac{n}{2}-1$. Using Lemmas \ref{L4} and \ref{L5}, we see that $ \alpha(G)=m_{L (G)} [0,1)\leq \gamma(G) \leq \alpha(G)$ which shows that $\gamma(G)=\alpha(G)=\frac{n}{2}$. Therefore, by Theorem 3 of {\cite{5R6}}, $G=HoK_1$ for some connected graph $H$. \qed

In the following theorem, we show that we can relax the condition $\alpha(G)=\frac{n}{2}$ in Theorem \ref{T5} for the class of bipartite graphs.

 \begin{theorem} \label{T6}  Let $G$ be a connected bipartite graph with $n$ vertices having independence number $\alpha(G)$ and diameter $d(G)\leq 2$. Then, $m_{D^{L} (G)} (2n-1,2n )= \alpha(G)-1$    if and only if $G=HoK_1$ for some connected graph $H$.
\end{theorem}
 \noindent {\bf Proof.} Assume that $G=HoK_1$, for some connected graph $H$. Then the proof follows by Theorem \ref{T5}. So let $m_{D^{L} (G)} (2n-1,2n )= \alpha(G)-1$. Using Theorem \ref{T5}, it is sufficient to show that $\alpha(G)=\frac{n}{2}$. If possible, let the two parts of $G$ have different orders. Then, using Lemmas \ref{L4} and  \ref{L5}, we have
 $$ \gamma(G)<\frac{n}{2}<\alpha(G)=m_{D^{L} (G)} (2n-1,2n )+1= m_{L (G)} [0,1)\leq \gamma(G),$$
 which is a contradiction. Therefore, the two parts of $G$ have the same order. Now, if $ \alpha(G)> \frac{n}{2}$, then by Lemma \ref{L7}, $m_{D^{L} (G)} (2n-1,2n )\leq \alpha(G)-2$, a contradiction. Hence $\alpha(G)\leq \frac{n}{2}$. Since the partite sets have the same order, we get $\alpha(G)=\frac{n}{2}$.\qed

 \noindent {\bf {Remark.}} From the above theorem, we see that if $G$ is a connected bipartite graph with $n$ vertices, having independence number $\alpha(G)$ and diameter $d\leq 2$ satisfying either of the conditions  (i)  $G=HoK_1$ for some connected graph $H$, or (ii) $m_{D^{L} (G)} (2n-1,2n )= \alpha(G)-1$, then $\alpha(G)=\frac{n}{2}$ and $n$ is even.\\

The following theorem shows that the number of distance Laplacian eigenvalues of the graph $G$ in the interval $[0,dn]$ is at least $d+1$.

\begin{theorem} \label{T10}  If $G$ is a connected graph of order $n$ having diameter $d$, then $$ m_{D^{L} (G )} \big[0,dn]\geq d+1.$$
\end{theorem}
\noindent {\bf Proof.} We consider the principal submatrix, say $M$, corresponding to the vertices $v_1 ,v_2 ,\dots, v_{d+1}$ which belong to the induced path $P_{d+1}$ in the distance Laplacian matrix of $G$. Clearly, the transmission of any vertex in the path $P_{d+1}$ is at most $\frac{d(2n-d-1)}{2}$, that is, $Tr(v_i )\leq \frac{d(2n-d-1)}{2}$, for all $i=1,2,\dots,d+1$. Also, the sum of the off diagonal elements of any row of $M$ is less than or equal to $\frac{d(d+1)}{2}$.  Using Lemma \ref{L2}, we conclude that the maximum eigenvalue of $M$ is at most $dn$. Using Fact 1 and Theorem \ref{T7}, there at least $d+1$ distance Laplacian eigenvalues of $G$ which are greater than or equal to $0$ and less than or equal to $dn$, that is,  $ m_{D^{L} (G )} \big[0,dn]\geq d+1.$ \qed

From Theorem \ref{T10},  we get the following observation after using Inequality (2.1).

\begin{corollary} \label{C4} Let $G$ be a connected graph of order $n$ having diameter $d$. If $dn<2Tr_{max}$, then $$ m_{D^{L} (G )} \big(dn,2Tr_{max}]\leq n- d-1.$$
\end{corollary}

\section{Concluding Remarks}

 In the entire generality, we believe it is hard to characterize all the graphs satisfying the bounds given in Theorem \ref{T1} and Theorem \ref{T8}. Also in Theorem \ref{T5}, we characterized graphs with diameter $d\leq 2$ satisfying $m_{D^{L} (G)} (2n-1,2n )= \alpha(G)-1=\frac{n}{2}-1$ and we left the case when $d\geq 3$.  So, the following problems will be interesting for future research.\\
 {\bf Problem 1.} {\it Determine the classes of graphs $\vartheta$ for which  $m_{D^{L} (G)} [n,n+\alpha(G))= n-\alpha(G)$, for any $G\in \vartheta$. } \\
 {\bf Problem 2.} {\it Determine the classes of graphs $\vartheta$ for which  $m_{D^{L} (G)} [n,n+2)= \chi-1$, for any $G\in \vartheta$. }\\
{\bf Problem 3.} {\it Determine the classes of graphs $\vartheta$ for which  $m_{D^{L} (G)} (2n-1,2n )= \alpha(G)-1=\frac{n}{2}-1$, for any $G\in \vartheta$ with $d\geq 3$. }\\

\noindent{\bf Data availibility} Data sharing is not applicable to this article as no data sets were generated or analyzed
during the current study.

\end{document}